\documentclass[11pt]{article}
\usepackage{latexsym,amsmath,amssymb}
\newtheorem{theorem}{Theorem}[section]

\newtheorem{corollary}{Corollary}[section]

\numberwithin{equation}{section}

\begin{document}

\title{Proposed theorems for lifts of the extended almost complex structures on the complex manifold}
\author{Mohammad Nazrul Islam Khan \\
Department of Computer Engineering, College of Computer \\
Qassim University, Buraydah, Saudi Arabia\\
Email: m.nazrul@qu.edu.sa, mnazrul@rediffmail.com\\}
\date{}
\maketitle
\begin{abstract}
It is well known that the tensor field $J$ of type (1,1) on the manifold $M$ is an almost complex structure if $J^2=-I, I$ is an identity tensor field and the manifold $M$ is called the complex manifold. Let $^kM$ be the $k$ order extended complex manifold of the manifold $M$. A tensor field $J_k$ on $^kM$ is called extended almost complex structure if $(J_k)^2=-I$. The present paper aims to study the higher order complete and vertical lifts of the extended almost complex structures on an extended complex manifold $^kM$. The proposed theorems on the Nijenhuis tensor of an extended almost complex structure $J_k$ on the extended complex manifold $^kM$ are proved. Also, a tensor field $\tilde{J_k}$ of type (1,1) is introduced and shows that it is an extended almost complex structure. Finally, the Lie derivative concerning higher-order lifts is studied and basic results on the almost analytic complex vector concerning an extended almost complex structure on $^kM$ are investigated.  
\end{abstract}




\section{Introduction}
The lifting theory has precious vicinity in the differentiable geometry due to the fact it is feasible to generalize to the geometric structures on any manifold to extensions through lift function. The geometric structures; for example, an almost complex structure, an almost product structure, on the base manifold admit lifts, namely, the complete and vertical lifts, to the canonical extended manifold. Tekkoyun et al \cite{TCG} studied the geometric structures of an extended complex manifold $^kM$ of $k$-th order of the complex manifold $M$ and established higher-order vertical and complete lifts of functions, vector fields, and 1-forms on $M$ to $^kM$. Tekkoyun and Civelek \cite{MT, TC} investigated higher order complete, vertical and horizontal lifts of the complex structures on the complex manifold $M$ to the extended complex manifold $^kM$. Das and the author \cite{DK} have studied almost $r$-contact structures on the tangent bundle using the complete and vertical lifts. The geonmtric structures such as  almost complex, metallic structure, almost Hermitian structure, etc. on the manifold to  the tangent bundle are studied by earlier investigators \cite{MGK, GM, ASG, SB, KCH, MJ, TO, AAS}. 

The major focus of this work can be summarized as follows:
\begin{itemize}
\item The proposed theorem proves that the higher order complete lift $F^{c^k}$ of an almost complex structure $F$ of type (1,1) is an extended almost complex structure on $^kM$. 
\item The Nijenhuis tensor of an extended almost complex structure $F^{c^k}$ on $^kM$ is determined.

\item A tensor field $\tilde{J_k}$ of type (1,1) on $^kM$ is introduced and shows that it is an extended almost complex structure on $^kM$.

\item The basic results on the almost analytic complex vector concerning an extended almost complex structure on $^kM$ are investigated.  
\end{itemize}

The paper is organized as follows: Section 2 describes a brief account of an extended almost complex manifold, higher-order lifts of complex function, vector field, 1-form, and tensor field of type (1,1). In Section 3, higher-order lifts of the Nijenhuis tensor and integrability conditions are calculated. In Section 4, a tensor field  
$\tilde{J_k}$ of type (1,1) is introduced on the extended complex manifold $^kM$ and shows that it is an extended almost complex structure on $^kM$. Finally, the Lie derivative concerning higher-order lifts is studied and some basic results on the almost analytic complex vector concerning an extended almost complex structure on $^kM$ are investigated.

\section{Preliminaries}
Let $M$ be $2m$-real dimensional manifold and $^kM$ its $k$-th order extended manifold. A tensor field $J_k$ of type (1,1) on $^kM$ is called an extended almost complex structure on $^kM$ if $J_k$ is endomorphism of the tangent space $Tp(^kM)$ such that $(J_k)^2 = -I$ at every point $p$ of $^kM$. An extended manifold $^kM$ with an extended almost complex structure $J_k$ is called an extended almost complex manifold \cite{TCG}. 

Suppose $(x^{ri}, y^{ri})$ be a system of real coordinate defined at any point $p$ of $^kM$ covered by neighborhood $^kU$ and $\{\frac{\partial}{\partial x^{ri}}, \frac{\partial}{\partial y^{ri}}\}$ be  the natural base over the real field $R$ of the tangent space $Tp(^kM)$ of $^kM$. The manifold $^kM$ is called an extended complex manifold if $^kM$ is covered by neighborhood $^kU$ and a local coordinate system $(x^{ri}, y^{ri})$ defined in the neighborhood $^kU$ such that
\begin{eqnarray}
J_k(\frac{\partial}{\partial x^{ri}})&=&\frac{\partial}{\partial y^{ri}},\\
J_k(\frac{\partial}{\partial y^{ri}})&=& - \frac{\partial}{\partial x^{ri}}.
\end{eqnarray}
If $k = 0, J_0$ is called an almost complex structure and a manifold $^0M=M$ with an almost complex structure $J_0$ is said to be the complex manifold.

Let $z^{ri} = x^{ri} + i y^{ri}, i = \sqrt{-1}$ be an extended complex local coordinate system on a neighborhood $^kU$ of any point $p$ of $^kM$. Then $\frac{\partial}{\partial z^{ri}} = \frac{1}{2}\{\frac{\partial}{\partial x^{ri}} - i\frac{\partial}{\partial y^{ri}}\} , \frac{\partial}{\partial \bar{z}^{ri}} = \frac{1}{2}\{\frac{\partial}{\partial x^{ri}} + i\frac{\partial}{\partial y^{ri}}\}$ and the endomorphism $J_k$ is given by
\begin{equation}
J_k(\frac{\partial}{\partial z^{ri}})=i\frac{\partial}{\partial z^{ri}},J_k(\frac{\partial}{\partial \bar{z}^{ri}})=-i\frac{\partial}{\partial \bar{z}^{ri}}.
\end{equation}

Let $M$ be any complex manifold and $^kM$ its $k$-th order extension. If $f$ is a function on $M$, then the functions $f^{v^k}$ and $f^{c^k}$ denote the  vertical and complete lifts of the function $f$ on $^kM$, respectively and given by \cite{TCG, OLP}   
\begin{equation}
f^{v^k} = f\circ \tau M\circ \tau ^2M \circ ... \circ \tau ^{k-1}M, 
\end{equation}
and
\begin{equation}
f^{c^k} = \dot{z}^{ri}(\frac{\partial f^{c^k}}{\partial z^{ri}})^v +\dot{\bar{z}}^{ri}(\frac{\partial f^{c^k}}{\partial \bar{z}^{ri}})^v,
\end{equation}
where $\tau ^{k-1}M\rightarrow^k M ^{k-1}M$ is a canonical projection. 

The extended properties for the vertical and complete lifts of the complex functions are given as follows:
\begin{eqnarray}
i)~~ (f + g)^{v^r} &=& f^{v^r} + g^{v^r} , (f.g)^{v^r}  = f^{v^r}.g^{v^r}\nonumber\\ 
ii)~~ (f + g)^{c^r} &=& f^{c^r} + g^{c^r}, (f.g)^{c^r} =\sum^r_{j=0}C^r_j f^{c^{r-j}v^j}.g^{c^j v^{r-j}},
\end{eqnarray}
where $f$ and $g$ are the complex functions and $C^r_j$ is the combination.

Let $X$ be a complex vector field with the local expression $X = Z^{0i}\frac{\partial }{\partial z^{ki}}+\bar{Z}^{0i}\frac{\partial }{\partial \bar{z}^{ki}}$. Then the local expression of the vertical and complete lifts of $X$ to $^kM$ are, respectively, given by     
\begin{equation}
X^{v^k} = (Z^{0i})^{v^k}\frac{\partial }{\partial z^{ki}}+(\bar{Z}^{0i})^{v^k}\frac{\partial }{\partial \bar{z}^{ki}}.
\end{equation}
and
\begin{equation}
X^{c^k} =C^r_j(Z^{0i})^{v^{k-r}c^r}\frac{\partial }{\partial z^{ri}}+C^r_j(\bar{Z}^{0i})^{v^{k-r}c^r}\frac{\partial }{\partial \bar{z}^{ri}}. 
\end{equation}

The extended properties for the vertical and complete lifts of the complex vector fields are given as follows:
\begin{eqnarray}\label{XYC}
i)~~ (X + Y)^{v^r} &=& X^{v^r} + Y^{v^r} , (X + Y)^{c^r} = X^{c^r} + Y^{c^r},\nonumber\\ 
ii)~~ X^{v^k}(f^{c^k}) &=& (Xf)^{v^k},~~~X^{c^k}(f^{c^k}) = (Xf)^{c^k},\nonumber\\
iii)~~~~ (fX)^{v^r} &=& f^{v^r}.X^{v^r}, (fX)^{c^r}  = \sum^r_{j=0}C^r_j f^{c^{r-j}v^j}.X^{c^j v^{r-j}},\\
iv)~~ X^{v^r}(f^{v^r})&=&0,~~X^{c^r}(f^{c^r})=(Xf)^{c^r},~X^{c^r}(f^{v^r})=X^{v^r}(f^{c^r})=(Xf)^{v^r},\nonumber\\
v)~~ [X^{v^r},Y^{v^r}]&=&0, [X^{c^r},Y^{c^r}]=[X,Y]^{c^r},[X^{v^r},Y^{c^r}]=[X^{c^r},Y^{v^r}]=[X,Y]^{v^r},\nonumber
\end{eqnarray}
where $X,Y$ are the complex vector fields and $f$ is the complex function.

Let $\alpha$ be a complex 1-form  with the local expression $\alpha = \alpha_{0i}dz^{0i} + \bar{\alpha}_{0i}\bar{dz}^{0i}$. Then the local expression of the vertical and complete lifts of $\alpha$ to $^kM$ are, respectively, given by     
\begin{equation}
\alpha^{v^k} = (\alpha_{0i})^{v^k}dz^{0i} + (\bar{\alpha}_{0i})^{v^k}\bar{dz}^{0i},
\end{equation}
and
\begin{equation}
\alpha^{c^k} = (\alpha_{0i})^{c^{k-r}v^r}dz^{ri} +  (\bar{\alpha}_{0i})^{c^{k-r}v^r}d\bar{z}^{ri}.
\end{equation}
The extended properties for the vertical and the complete lifts of complex 1-forms are given as follows \cite{TCG}:
\begin{eqnarray}
i) (\alpha+\lambda)^{v^r}&=&\alpha^{v^r}+\lambda^{v^r}, (\alpha+\lambda)^{c^r}=\alpha^{c^r}+\lambda^{c^r},\nonumber\\
ii)~~~ (f\alpha)^{v^r}&=&f^{v^r}\alpha^{v^r},(f\alpha)^{c^r}=\sum^r_{j=0}C^r_jf^{c^{r-j}v^j}\alpha^{c^rv^{r-j}}.  
\end{eqnarray}

Let $M$ be any complex manifold and $^kM$ its $k$-th order extension. If $F$ be a tensor field of type (1,1). Then
\begin{eqnarray}
\alpha^{c^k}(F^{c^k})&=&(\alpha F)^{c^k},~~\alpha^{v^k}(F^{v^k})=(\alpha F)^{v^k},\nonumber\\
F^{c^k}(X^{c^k})&=&(FX)^{c^k},~~F^{v^k}(X^{c^k})=(FX)^{v^k},
\end{eqnarray}
where $X$ and $\alpha$ are a vector field and a 1-form respectively.

\section{Some calculations on the Nijenhuis tensor on the extended complex manifold}
In this section, the basic properties of the tensor field of type (1,1) on an extended complex manifold $^kM$ are studied. The proposed theorem proves that  the higher order complete lift $F^{c^k}$ of a tensor field $F$ of type (1,1) is an extended almost complex structure to an extended complex manifold $^kM$ if $F$ is an almost complex structure on $M$. The Nijenhuis tensor of an extended almost complex structure $F^{c^k}$ on the extended complex manifold $^kM$ is determined.
 
Let $F$ and $G$ be the complex tensor fields of type (1,1) and consider them as fields of linear endomorphisms of tangent spaces of complex manifold $M$. Let the eteration of the endomorphisms $F$ and $G$ be denoted by $FG$ i.e. $(FG)X=F(GX)$, where $X$ being an arbitrary complex vector field \cite{YI}.
\begin{theorem} Let $M$ be any complex manifold and $^kM$ its $k$-th order extension. If $F$ and $G$ be the complex tensor fields of type (1,1). Then 
\begin{eqnarray}
(a)~~~(FG)^{c^k}=F^{{c^k}G^{c^k}},\nonumber\\
(b)~~~~~I^{c^k}=I,
\end{eqnarray}
where $I$ is an identity tensor field.
\end{theorem}
\textit{Proof.}
 (i) Let $M$ be any complex manifold and let $F$ and $G$ be complex tensor fields of type (1,1). Then
\begin{eqnarray}
(FG)^{c^k}X^{c^k}&=&((FG)X)^{c^k}=(F(GX))^{c^k}=F^{c^k}(GX)^{c^k}\nonumber\\
&=&F^{c^k}G^{c^k}X^{c^k}=(F^{c^k}G^{c^k})X^{c^k}\nonumber\\
\Rightarrow(FG)^{c^k}&=&F^{c^k}G^{c^k},
\end{eqnarray}
for an arbitrary complex vector field $X$.

Second part of the theorem is proved by using similar devices.

Hence the theorem is proved.\\

\textbf{Definition 1}
Let $M$ be any complex manifold and $^kM$ its $k$-th order extension. 
If the tensor field $F$ of type (1,1) on $M$ satisfies the equation 
\begin{equation}\label{FI}
F^2+I=0, 
\end{equation} 
where $I$ is an identity tensor field. Then $F$ is called an almost complex structure on $M$ \cite{KYI}.

\begin{theorem}
Let $F$ be a tensor field  of type (1,1) on $M$. Then $F^{c^k}$ defines an extended almost complex structure to $^kM$ if and only if $F$ defines an almost complex structure on $M$.
\end{theorem}
\textit{Proof.} Let $F$ be an almost complex structure on $M$, then  
$$F^2+I=0.$$
Taking the complete lifts on both sides the above equation to $^kM$, then
\begin{eqnarray*}
(F^2+I)^{c^k}&=&0,\\
(F^2)^{c^k}+I&=&0,\\
(F^2)^{c^k}&=&-I\Rightarrow (F^{c^k})^2=-I.
\end{eqnarray*}
Therefore, $F^{c^k}$ defines an extended almost complex structure to $^kM$.

Hence the theorem is proved.\\\\

The Nijenhuis tensor $N_F$ of $F$ is a tensor field of type (1,2) defined by \cite{MNC, YI}
\begin{equation}\label{NF}
N_F(X,Y)=[FX,FY]+F^2[X,Y]-F[FX,Y]-F[X,FY],
\end{equation}
where $X$ and $Y$ are the complex vector fields in $M$.
\begin{theorem} Let $M$ be any complex manifold and $^kM$ its $k$-th order extension. Then 
$$(N_F)^{c^k}=N_{F^{c^k}},$$
where $X,Y$ are arbitrary complex vector fields.
\end{theorem}
\textit{Proof.} Taking lifts of both sides of equation (\ref{NF}) and using equation (\ref{XYC}), then
\begin{eqnarray*}
(N_F)^{c^k}&=&([FX,FY]+F^2[X,Y]-F[FX,Y]-F[X,FY])^{c^k},\\
&=&[F^{c^k}X^{c^k},F^{c^k}Y^{c^k}]+(F^{c^k})^2[X^{c^k},Y^{c^k}]-F^{c^k}[F^{c^k}X^{c^k},Y^{c^k}]-F^{c^k}[X^{c^k},F^{c^k}Y^{c^k}],\\
&=&N_{F^{c^k}}.
\end{eqnarray*}

Hence the theorem is proved.

\begin{theorem}
Let $F$ be a tensor field  of type (1,1) on $M$. Then an extended almost complex structure $F^{c^k}$ is integrable to $^kM$ if and only if an almost complex structure $F$ is integrable on $M$.  
\end{theorem}
\textit{Proof.} Let $N_{F^{c^k}}$ and $N_F$ be the Nijenhuis tensors of $F^{c^k}$ and $F$, respectively. Then 
\begin{equation}\label{NFC} 
N_{F^{c^k}}=(N_F)^{c^k},
\end{equation} 
since $F$ is integrable, then 
$$N_F=0.$$
From equation (\ref{NFC}), the obtained equation is 
$$N_{F^{c^k}}=0.$$
Thus, $F^{c^k}$ is integrable to $^kM$ if and only if $F$ is integrable on $M$.

Hence the theorem is proved.

\section{Extended almost complex structures}
In this section, a tensor field $\tilde{J_k}$ of type (1,1) on $^kM$ is introduced and shows that it is an extended almost complex structure on $^kM$.

Suppose that there is given, on a complex manifold $M$, a tensor field $F$ of type (1,1), a vector field $U$ and a 1-form $\omega$ satisfying 
\begin{eqnarray}\label{FU}
F^2&=&-I+U\otimes\omega,~~~ FU=0,\nonumber\\
\omega\circ F&=&0,~~~\omega(U)=0.
\end{eqnarray}
The structure $(F,U,\omega)$ is said to be an almost contact structure on $M$ \cite{YI}.
\begin{theorem} 
Let $M$ be any complex manifold and $^kM$ its $k$-th order extension  endowed with an almost contact structure $(F, U, \omega)$, then 
\begin{equation}
\tilde{J_k}=F^{c^k}+U^{v^k}\otimes\omega^{v^k}-U^{c^k}\otimes\omega^{c^k}
\end{equation}
 is an extended almost complex structure on $^kM$.
\end{theorem}
\textit{Proof.}
Taking complete lifts on both sides of equation (\ref{FU}) to $^kM$, then \cite{MN}
\begin{eqnarray}\label{FUC}
(F^{c^k})^2&=&-I+U^{v^k}\otimes\omega^{c^k}+U^{c^k}\otimes\omega^{v^k},\nonumber\\
F^{c^k}U^{v^k}&=&0,~~F^{c^k}U^{c^k}=0,~~\omega^{v^k}\circ F^{c^k}=0,~~\omega^{c^k}\circ F^{c^k}=0,\\
\omega^{v^k}(U^{v^k})&=&0,~~\omega^{c^k}(U^{c^k})=0,~~\omega^{c^k}(U^{v^k})=1,~~\omega^{v^k}(U^{v^k})=0.\nonumber
\end{eqnarray}
Let us introduce a tensor field $\tilde{J_k}$ of type (1,1) on $^kM$ as
\begin{equation}\label{JF}
\tilde{J_k}=F^{c^k}+U^{v^k}\otimes\omega^{v^k}-U^{c^k}\otimes\omega^{c^k}.
\end{equation}
Using equations (\ref{FUC}) and (\ref{JF}), the obtained equation is 
$$\tilde{J_k}^2=-I.$$

This shows that $\tilde{J_k}$ is an extended almost complex structure on $^kM$.

\begin{corollary} 
Let $M$ be a complex manifold and $\tilde{J_k}$ defined in equation (\ref{JF}) is an extended almost complex structure on $^kM$, then
\begin{eqnarray}\label{JXV}
\tilde{J_k}X^{v^k}&=&(FX)^{v^k}-(\omega(X))^{v^k}U^{c^k},\nonumber\\
\tilde{J_k}X^{c^k}&=&(FX)^{c^k}+(\omega(X))^{v^k}U^{v^k}-(\omega(X))^{v^k}U^{c^k},
\end{eqnarray}
where $X$ is a complex vector field and $\omega$ is a 1-form.
\end{corollary}
\begin{corollary}
Let $M$ be a complex manifold and $\tilde{J_k}$ is an extended almost complex structure on $^kM$ such that $\omega(X)=0$, then  
\begin{eqnarray}\label{JXP}
\tilde{J_k}X^{v^k}&=&(FX)^{v^k},~~~\tilde{J_k}X^{c^k}=(FX)^{c^k},\nonumber\\
\tilde{J_k}U^{v^k}&=&-U^{c^k},~~~\tilde{J_k}U^{c^k}=-U^{v^k},
\end{eqnarray}
where $X$ is a complex vector field and $\omega$ is a 1-form.
\end{corollary}
\section{Lie derivative with respect to higher order lifts on the extended complex manifold}
In this section, a study is done on the Lie derivative with respect to a complex vector field on $^kM$. The basic results on the almost analytic complex vector concerning an extended almost complex structure on $^kM$ are found.  

Let $M$ be any complex manifold and $X$ be a complex vector field. The differential transformation $\pounds_X$ is called Lie derivative with respect to $X$ if
\begin{equation}\label{LX}
\pounds_X f=Xf,~~~\pounds_X Y=[X,Y],
\end{equation}
where $f$ is complex function on $M$.	
	 	
The Lie derivative $\pounds_X F$ of a tensor field $F$ of type (1, 1) on $M$ with respect to a complex vector field $X$ is defined by \cite{YI}
\begin{equation}\label{LXF}
(\pounds_X F)Y=[X,FY]-F[X,Y],
\end{equation}	
where $[,]$ is Lie bracket.

\begin{theorem}
Let $M$ be the complex manifold and $^kM$ its $k$-th order extension. The Lie derivation $\pounds_X$ with respect to a complex vector field $X$, then           
\begin{align}\label{LXV}
\mbox{(a)}&~~~~ \pounds_{X^{v^r}} f^{v^r}=0,~~\pounds_{X^{v^r}} f^{c^r}=(\pounds_X f)^{v^r},\nonumber\\
\mbox{(b)}&~~~~\pounds_{X^{c^r}} f^{v^r}=(\pounds_X f)^{v^r},~~\pounds_{X^{c^r}} f^{c^r}=(\pounds_X f)^{c^r},\\
\mbox{(c)}&~~~~\pounds_{X^{v^r}} Y^{v^r}=0,~~\pounds_{X^{v^r}} Y^C=(\pounds_XY)^{v^r},\nonumber\\
\mbox{(d)}&~~~~\pounds_{X^{c^r}} Y^{v^r}=(\pounds_X Y)^{v^r},~~\pounds_{X^{c^r}} Y^{c^r}=(\pounds_X Y)^{c^r},\nonumber
\end{align}
where $f$ is a complex function in $M$.
\end{theorem}
\textit{Proof.} The proof follows in an obvious manner on using equation (\ref{XYC}).
\begin{theorem}
Let $M$ be the complex manifold and $^kM$ its $k$-th order extension. Let an extended almost complex structure $\tilde{J}_k$ defined by equation (\ref{JF}) and $\omega(Y)=0$ on $^kM$. The Lie derivation $\pounds_X$ with respect to a complex vector field $X$, then
\begin{align}\label{LXJ}
\mbox{(a)}&~~~~ (\pounds_{X^{v^r}} \tilde{J_k})Y^{v^r}=0,\nonumber\\
\mbox{(b)}&~~~~(\pounds_{X^{v^r}} \tilde{J_k})Y^{c^r}=\left((\pounds_XF)Y\right)^{v^r}-((\pounds_X\omega)Y)^{v^r}U^{c^r},\nonumber\\
\mbox{(c)}&~~~~(\pounds_{X^{c^r}} \tilde{J_k})Y^{v^r}=\left((\pounds_XF)Y\right)^{v^r}-((\pounds_X\omega)Y)^{v^r}U^{c^r,}\\
\mbox{(d)}&~~~~(\pounds_{X^{c^r}} \tilde{J_k})Y^{c^r}=\left((\pounds_XF)Y\right)^{c^r}+((\pounds_X\omega)(Y))^{v^r}U^{v^r}-((\pounds_X\omega)(Y))^{c^r}U^{c^r},\nonumber
\end{align}
where $F$ is a tensor field of type (1,1).
\end{theorem}
\textit{Proof.} Let $M$ be any complex manifold and $\omega(Y)=0$. Then 
\begin{eqnarray}
\mbox{(5.5a)}~~(\pounds_{X^{v^r}}\tilde{J_k})Y^{v^r}&=&\pounds_{X^{v^r}}\left[F^{c^k}+U^{v^k}\otimes\omega^{v^k}-U^{c^k}\otimes\omega^{c^k}\right]Y^{v^k},\nonumber\\
&-&\left[F^{c^k}+U^{v^k}\otimes\omega^{v^k}-U^{c^k}\otimes\omega^{c^k}\right]\pounds_{X^{v^r}}Y^{v^k}\nonumber\\
&=&\pounds_{X^{v^r}}F^{c^k}Y^{v^k}+\pounds_{X^{v^r}}(U^{v^k}\otimes\omega^{v^k})Y^{v^k}-\pounds_{X^{v^r}}(U^{c^k}\otimes\omega^{c^k})Y^{v^k},\nonumber\\
&=&\pounds_{X^{v^r}}(FY)^{v^k}+\pounds_{X^{v^r}}(\otimes\omega(Y))^{v^k}U^{v^k}-\pounds_{X^{v^r}}(\otimes\omega(Y))^{v^k}U^{c^k},\nonumber\\
&=&0, ~~~as ~\omega(Y)=0,~as \pounds_{X^{v^r}}Y^{v^k}=0.\nonumber\\
\mbox{(5.5b)}~~(\pounds_{X^{v^r}}\tilde{J_k})Y^{c^r}&=&\pounds_{X^{v^r}}\left[F^{c^k}+U^{v^k}\otimes\omega^{v^k}-U^{c^k}\otimes\omega^{c^k}\right]Y^{c^k}\nonumber\\
&-&\left[F^{c^k}+U^{v^k}\otimes\omega^{v^k}-U^{c^k}\otimes\omega^{c^k}\right]\pounds_{X^{v^r}}Y^{c^k},\nonumber\\
&=&\pounds_{X^{v^r}}F^{c^k}Y^{c^k}+\pounds_{X^{v^r}}(\otimes\omega^{v^k}(Y)^{c^k})U^{v^k}-\pounds_{X^{v^r}}(\otimes\omega^{c^k}(Y)^{c^k})U^{c^k},\nonumber\\
&-&F^{c^k}\pounds_{X^{v^r}}Y^{c^k}-\omega^{v^k}(\pounds_{X^{v^r}}Y^{c^k})U^{v^k}+\omega^{c^k}(\pounds_{X^{v^r}}Y^{c^k})U^{c^k},\nonumber\\
&=&((\pounds_XF)Y)^{v^r}-((\pounds_X\omega)Y)^{v^k}U^c.\nonumber
\end{eqnarray}
The other formulas the theorem are proved by using similar devices.
\begin{corollary}
If we put $Y=U$ , an extended almost complex structure $\tilde{J}_k$ defined by equation (\ref{JF}) and $\omega(Y)=0$, we have
\begin{align}\label{LXU}
\mbox{(a)}&~~~~ (\pounds_{X^{v^r}} \tilde{J_k})U^{v^r}=-(\pounds_XU)^{v^r},\nonumber\\
\mbox{(b)}&~~~~(\pounds_{X^{v^r}} \tilde{J_k})U^{c^r}=((\pounds_X F)U)^{v^r}-((\pounds_X \omega)(U))^{v^r}U^{c^r},\nonumber\\
\mbox{(c)}&~~~~(\pounds_{X^{c^r}} \tilde{J_k})U^{v^r}=((\pounds_X F)U)^{c^r}-[X,U]^{c^r}-((\pounds_X \omega)(U))^{c^r}U^{c^r},\\
\mbox{(d)}&~~~~(\pounds_{X^{c^r}} \tilde{J_k})U^{c^r}=((\pounds_X F)U)^{c^r}+[X,U]^{v^r}+((\pounds_X \omega)(U))^{v^r}U^{v^r}-((\pounds_X \omega)(U))^{c^r}U^{c^r},\nonumber
\end{align}
where $F$ is a tensor field of type (1,1).
\end{corollary}

\begin{theorem}
Let $M$ be any complex manifold and $^kM$ its $k$-th order extension endowed with an almost contact structure $(F, U, \omega)$. Then the vertical lift $X^{v^r}$ of a complex vector field $X$ given in $M$ is an almost analytic with respect to the extended almost complex structure $\tilde{J_k}$ defined by equation (\ref{JF}) on $^kM$ if and only if the conditions
 $$\pounds_XF=0,~~~~  \pounds_XU=0,~~~~   \pounds_X\omega=0$$
are satisfied on $M$. 
\end{theorem}
\textit{Proof.}
As a consequence of equations (\ref{LXJ}a), (\ref{LXJ}b), (\ref{LXU}a) and (\ref{LXU}b), it follows that the condition 
 		  $$\pounds_{X^{v^r}}J_k=0$$
is equivalent to the conditions
 $$\pounds_XF=0,~~~~  \pounds_XU=0,~~~~   \pounds_X\omega=0.$$

\begin{theorem}
Let $M$ be any complex manifold and $^kM$ its $k$-th order extension  endowed with an almost contact structure $(F, U, \omega)$. Then the vertical lift $X^{c^r}$ of a complex vector field $X$ given in $M$ is almost analytic with respect to the an extended almost complex structure $\tilde{J_k}$ defined by equation (\ref{JF}) on $^kM$ if and only if the conditions
 $$\pounds_XF=0,~~~~  \pounds_XU=CU, ~~~~  \pounds_X\omega=-C\omega,$$
$C$ being a non-zero constant, are satisfied on $M$.  
\end{theorem}
\textit{Proof.}
As a consequence of equations (\ref{LXJ}c), (\ref{LXJ}d), (\ref{LXU}c) and (\ref{LXU}d), it follows that the condition 
		 $$\pounds_{X^{c^r}}J_k=0$$
is equivalent to the conditions
$$\pounds_XF=0,~~~~  \pounds_XU=CU,~~~~   \pounds_X\omega=-C\omega.$$

\section*{\textsc{Conflict of  Interests}}
The authors declare that there is no conflict of interests.

\end{document}